\documentclass{article}

\usepackage{CJK}
\usepackage{amsmath}
\usepackage{amssymb}
\usepackage{amsfonts}
\usepackage{latexsym}
\usepackage{mathrsfs}

\CJKtilde                     

\begin{document}

\begin{CJK*}{GBK}{song}

\begin{center}
{\bf \large Partial regularities and $a^*$-invariants of Borel type ideals}

\vspace{4mm} Dancheng Lu and Lizhong Chu

\vspace{3mm}

  {\footnotesize \it Department of Mathematics, Soochow University,  Suzhou 215006, P.R. China}

\vskip 30pt

\begin{minipage}{14 true cm}

{\small\noindent{\bf Abstract:}\ We express the Partial regularities and $a^*$-invariants of a Borel type ideal in terms of its irredundant irreducible decomposition. In addition we consider the
behaviours of those invariants under intersections and sums. \vskip
0.1cm

{\bf Keywords:}\  Borel type ideals, regularity, stable \vskip 1mm

{\bf 2010 MR Subject Classification}:\  Primary 13D02,  Secondary 13A02}

\end{minipage}
\end{center}

\begin{figure}[b]
\vspace{-4mm}
\rule[-2.5truemm]{5cm}{0.1truemm}\\[2mm]
{\small

\ \  \ \ This work is supported by Natural Science Foundation of
China (11201326).}

\end{figure}

\vskip 4pt

\vspace{3mm}\section{ Introduction}

Let $k$ be an infinite field and $S=k[x_1,\cdots,x_n]$ be a
polynomial ring over $k$. Herzog, Popescu and Vladoiu \cite{HPV}
define a monomial ideal $I$ in $S$ to be of {\it Borel type} if it
satisfies  the following property:
$$ I:(x_1,\cdots,x_i)^{\infty}=I:x_i^{\infty}, \forall i=1,\cdots,n.  $$
We mention that this concept appears also in \cite[Definition
1.3]{CS} as   weakly  stable ideal and in
\cite[Definition 3.1]{BG} as  nested type ideal. The class
of ideals is rather large: It includes Borel fixed ideals (i.e.,
strongly stable ideals if $\mathrm{char}(k)=0$ or $p$-Borel ideals
if $\mathrm{char}(k)=p$. See \cite{BS}) and stable ideals. Many
authors are interested in the regularity of those ideals (see e.g.
\cite{HPV,BG,C,C1}.) For instance, if we call
$sat(I)=\max\{j:H_m^0(R/I)_j\neq 0\}$ to be the satiety of an ideal
$I$,  it is proved in \cite{HPV} and \cite{BG} independently   that
the regularity of an ideal of Borel type is  the maximum of satieties of
of some related graded ideals. Moreover it is proved in \cite{HPV}
that every monomial ideal of Borel type is sequentially Cohen-Macaulay.

\vspace{2mm} Let $M$ be a finitely generated graded $S$-module and $\mathfrak{m}$ be an the maximal graded ideal of $S$.
Note that $H_{\mathfrak{m}}^i(M)$ is an artinian graded module. As usual, we set
$$a_i(M)=\max\{j\, |\, H_{\mathfrak{m}}^i(M)_j\neq
0\}=\max\{j\, |\, \mathrm{Ext}^{n-i}_S(M,S)_{-n-j}\neq 0\}.$$ The
following invariants (called partial regularities and $a^*$-invariants)  are studied in \cite{T}.
$$\mathrm{reg}_t(M)=\max\{a_i(M)+i: i\leq t\}.$$
$$a^*_t(M)=\max\{a_i(M):i\leq t\}.$$

These invariants carry important information on the structure of
$M$. For instance, it is proved in \cite[Theorem 3.1]{T} that
$\mathrm{reg}_t(M)=\max\{b_i-i:i\geq n-t\}$ and
$a_t^*(M)=\max\{b_i:i\geq n-t\}-n$, where $b_i$ is the maximal
degree of generators in $F_i$ and where $0\rightarrow
F_s\rightarrow\cdots \rightarrow F_1\rightarrow F_0\rightarrow
M\rightarrow 0$ is the minimal free resolution of $M$. We see that
the usual (Castelnuovo-Mumford) regularity and $a^*$-invariant are special cases of these
invariants.

\vspace{2mm} In this short note we will study those partial
invariants a monomial ideal of Borel type.  Our main contribution is
Proposition 2.5 and Theorem 2.6, which say that all $a_i(S/I),\, (i=0, 1, \cdots, n)$ can be expressed explicitly in terms of its irredandunt irreducible decomposition if $I$ is of Borel type.
We see that
\cite[Corolary 3.7]{BG} follows immediately from these results. In
addition, we recover some results which appear in \cite{C,C1}.

\vspace{2mm}
 Let $I$
be a monomial ideal of $S$ and $u$ a monomial. We set $$m(u)=\max\{i|x_i \mbox{ divides }  u\}$$ and
$$m(I)=\max\{m(u)|u\in G(I)\},$$ where $G(I)$ is the set of minimal
generators of $I$. Let $N$ is a finitely generated graded
$S$-module. We use $s(N)$ to denote $\max\{i\, | \, N_i\neq 0\}$, the largest non-vanishing degree of $N$.

\vspace{2mm}

We conclude this section by the following characterizations of a monomial ideal of
Borel type.

\vspace{2mm} {\bf \noindent Theorem  1.1.} {\it Let $I$ be a
monomial ideal with $d=\dim (S/I)$. Then the following statements
are equivalent:}

(1) {\it $I$ is of Borel type.}

(2) {\it $x_n,\cdots,x_{1}$ is an almost regular sequence  on $S/I$.}

(3) {\it $H_1(x_n,\cdots,x_{i},S/I)$ has the finite length for $i=1,\cdots,n$. }

(4) {\it  $I_{\geq e}$ is stable for all $e>>0$, where $I_{\geq e}$ is the ideal generated by monomials of degree $\geq e$ from $I$.}

 \vspace{2mm}
 \noindent{\bf Proof.}
$\mathrm{ (1)} \Leftrightarrow \mathrm{(2)}$:  It follows from \cite[Theorem 3.1]{BG}.

$\mathrm{ (2)} \Leftrightarrow \mathrm{(3)}$:  It follows from \cite[Proposition 4.3.5]{HH}.

$(3)\Rightarrow (4)$: By the condition (3), there exists a positive integer $t$ such that $[H_1(x_n,\cdots,x_i;R/I)])_j=0$ for any  $1\leq i\leq n$ and $j\geq t$.  Then, for all $e\geq t$ and $i=1,\cdots,n$, we have that $$H_1(x_n,\cdots,x_i; S/I_{\geq e})=H_2(x_n,\cdots,x_i;I_{\geq e})$$
 $$=[H_2(x_n,\cdots,x_i;I)]_{\geq e}=[H_1(x_n,\cdots,x_i;S/I)]_{\geq e}=0.$$ It follows that $I_{\geq e}$ is stable by \cite[Theorem 7.2.2]{HH}.

$(4)\Rightarrow (3)$:
 If (3) does not hold, then there exists  $1\leq i\leq n$ such that $$\mathrm{dim }(H_1(x_n,\cdots,x_i; S/I))>0.$$ It follows that $\mathrm{dim }(H_1(x_n,\cdots,x_i; R/I_{\geq e}))>0$ for all $e>0$ and thus $I_{\geq e}$ is not  stable for any $e>0$ by \cite[Theorem 7.2.2]{HH}. $\Box$

 \vspace{2mm}

We remark that the equivalence between (1) and (4) has been proved in [9].

\vspace{4mm}\section{Main Results}

Let $I$ be a Borel type monomial ideal in $S$.
 Recall that the sequential chain $I=I_0 \subseteq I_1\subseteq\cdots\subseteq I_r=S$ of $I$ is
 defined in \cite{HPV} recursively:  $I_0=I$. Suppose $I_l$ is already defined.
 If $I_l=S$ then the chain ends. Otherwise $I_{l+1}=(I_l: x_{n_l}^{\infty})$, where $n_l=m(I_l)$.

 Notice that $r\leq n$ since $n\geq n_0>\cdots n_l>n_{l+1}>\cdots>n_{r-1}$.
 Let $J_l$ be the ideal generated by $G(I_l)$ in $k[x_1,\cdots,x_{n_l}]$ for $l=0,1,\cdots,r-1$. Then $I_{l+1}/I_l\cong J_l^{sat}/J_l[x_{n_l+1},\cdots,x_n]$ is a Cohen-Macaulay $S$-module of dimension $n-n_l$.
 This shows $S/I$ is a sequentially Cohen-Macauly module. Moreover we have:

 \vspace{2mm}
{\bf \noindent Proposition  2.1.} (\cite[Corollary 2.7.]{HPV}) {\it
Let $I$ be a monomial ideal of Borel type. With notations
$I_i's,J_i's$ as above we have $a_{n-j}(S/I)=-\infty$ if $j\notin
\{n_0,\cdots,n_{r-1}\}$ and
$$a_{n-n_i}(S/I)=s(J_{i}^{sat}/J_{i})-n+n_i, \ \forall \ i=0,1,\cdots, r-1.$$  }

\vspace{2mm}

{\bf \noindent Remark 2.2.} We use the convention that
$\max\{\emptyset\}=-\infty$, and  $s(N)=-\infty$ if $N=0$. Thus
$a_{n-j}(S/I)=-\infty$, rather than $0$, which appeared in the original
paper \cite{HPV}  if $j\notin \{n_0,\cdots,n_{r-1}\}$.

\vspace{2mm}
 By \cite[Corollary 2.4]{BG} if I is a monomial ideal satisfying
 $I:(x_1,\cdots,x_n)=I:x_n$, then $s(I^{sat}/I)$ is the maximal degree of
 minimal generators which involve  $x_n$. By this result and
 Proposition 2.1, we recover a result of a strongly stable monomial ideal given in \cite{T}, where a monomial ideal $I$ is called strongly stable if $I:(x_1,\cdots,x_i)=I:x_i$, $\forall i=1,2, \cdots, n$.

\vspace{2mm}\noindent {\bf  Corollary  2.3.} (\cite[Theorem 2.4]{T})
{\it Let $I$ be a strongly stable monomial ideal. Then}

(1) {\it $\mathrm{reg}_t(I)$ is the maximal degree of minimal
generators $u$ with $u\in G(I)$ and  $m(u)\geq n-t$.}

(2) {\it $a^*_t(I)$ is the maximal number of $deg(u)-n+m(u)$, where $u$ is a
minimal generator  of $I$ with $m(u)>n-t.$}

\vspace{2mm}
 \noindent{\bf Proof.} Let $\delta_i(I)$ denote the maximal degree
 of minimal generators
  $u$ of $I$ with $m(u)=i$. By the
 construction of $J_l's$ we see that
 $\delta_{n_0}(J_0)=\delta_{n_0}(I)$ and  $$\delta_{n_l}(I)\leq \delta_{n_l}(J_l)\leq \max\{\delta_{n_0}(I),\cdots
 \delta_{n_l}(I)\},\  \forall \, l=1,\cdots,r-1.   \qquad (\ast)$$

 Fix $t\in \{1,\cdots,n\}$, there exists $l$ such that $n-n_{l-1}<
 t-1\leq n-n_{l}$. Since $a_0(I)=-\infty$, $a_i(I)=a_{i-1}(S/I), i=1,\cdots,n-1$ and
  $\max\{-n,a_n(I)\}=\max\{-n,a_{n-1}(S/I)\}$,
 it follows that $$a^*_t(I)=\max\{a_0(I),\cdots,a_t(I)\}=\max\{a_0(S/I),\cdots, a_{t-1}(S/I)\}$$
 $$=\max\{a_{n-n_0}(S/I),\cdots,a_{n-n_l}(S/I)\}$$
 $$=\max\{\delta_{n_0}(J_0)-n+n_0,\cdots,
 \delta_{n_l}(J_l)-n+n_l\},$$ where the last inequality follows
 from \cite[Corollary 2.4]{BG}. In view of $(\ast)$ we have $$a^*_t(S/I)=\max\{\delta_{n_0}(I)-n+n_0,\cdots,
 \delta_{n_l}(I)-n+n_l\}.$$
 Note that $\delta_l(I)=-\infty$ if $l\notin
 \{n_0,\cdots,n_{r-1}\}$, we prove (2). (1) can be proved
 similarly. $\Box$

 \vspace{2mm} It is well-known that every monomial ideal $I$ has a unique irredundant decomposition  $I=q_1\cap q_2\cap \cdots \cap q_r$ , where each $q_i$ is an irreducible monomial ideal, i.e., an ideal generated by powers of variables (See \cite[Theorem 5.1.17]{MS}). For the writing of an irreducible monomial ideal we adopt a notation from \cite{MS}: If ${\bf{b}}=(b_1,b_2,\cdots, b_n)\in {\mathbb{N}}^n$, then ${\mathfrak{m}}^{\bf{b}}$ denotes the irreducible monomial ideal $(x_i^{b_i}: b_i\geq 1).$  Set ${|{\bf{b}}|}=b_1+b_2+\cdots +b_n$.

 \vspace{2mm}

  {\bf \noindent
Lemma 2.4.} {\it Let $I={\mathfrak{m}}^{\bf{b}}$ with
$b_i>0,\ i=1,\cdots,n$. Then $s(I^{sat}/I)={|{\bf{b}}|}-n$.}

\vspace{2mm}
 \noindent{\bf Proof.} Straightforward check.

 \vspace{2mm}
{\bf \noindent Proposition  2.5.} {\it Let $I$ be a monomial ideal
of Borel type and let $I=q_1\cap q_2\cap \cdots \cap q_r$ be a
 irredundant irreducible decomposition of $I$.
Then
$$s(I^{sat}/I)=\max\{s(q_1^{sat}/q_1),\cdots, s(q_r^{sat}/q_r)\}.$$}
 \noindent{\bf Proof.}
  Case 1. If $(x_1,\cdots,x_n)\notin \mathrm{Ass}(I)$, then $I^{sat}=I:x_n^{\infty}=I$ and $q_i^{sat}=q_i,i=1,\cdots, r$. Hence $s(I^{sat}/I)=s(q_i^{sat}/q_i)=-\infty, 1\leq i\leq r,$ as required.

  Case 2. If $(x_1,\cdots,x_n)\in \mathrm{Ass}(I)$, we let $$K=\cap \{q_i|\sqrt{q_i}=(x_1,\cdots,x_n),1\leq i\leq r\}$$ and $$L=\cap \{q_i|\sqrt{q_i}\neq (x_1,\cdots,x_n),1\leq i\leq r\}.$$  By permutation of subscripts we can  assume $K=q_1\cap\cdots \cap q_s$ and $L=q_{s+1}\cap \cdots\cap q_r$ where $s\leq r$. Since $K^{sat}=S$, we have $s(K^{sat}/K)=\max\{t|K_t\neq S_t\}=\max\{t|(q_i)_t\neq S_t,\forall 1\leq i\leq s\}=\max\{s(q_i^{sat}/q_i)|1\leq 1\leq s\}.$ Denote $T=s(K^{sat}/K)$. We only need to prove that $s(I^{sat}/I)=T$ since $s(q_i^{sat}/q_i)=-\infty$ for each $i\in \{s+1,\cdots, r\}.$

Note that $s(I^{sat}/I)=s(L/K\cap L)=\max\{t|L_t\nsubseteqq K_t\}.$  If $t>T$, then $L_t\subseteq K_t=S_t$. It follows that $s(I^{sat}/I)\leq T$. Thus the results follows if one shows that $L_T\nsubseteqq K_T$. By Lemma 2.4  there exists $1\leq i\leq s$, say $i=s$, such that $q_s=(x_1^{b_1},\cdots,x_n^{b_n})$ and $T=b_1+\cdots+b_n-n$ where $b_i>0, i=1,\cdots,n.$ For each $s+1\leq j\leq r$, since $q_{j}\nsubseteqq q_s$ there is $x_{k_j}^{a_j}\in G(q_j)$ such that $x_{k_j}^{a_j}$ divides $x_1^{b_1-1}\cdots x_n^{b_n-1}$. Let $u$ be the least common multiple of $x_{k_j}^{a_j},j=s+1,\cdots,r$. Since $u\in L$ we have $x_1^{b_1-1}\cdots x_n^{b_n-1}\in L_T$ and so $L_T\nsubseteqq K_T$, completing the proof. $\Box$

\vspace{2mm}

By virtue of Propositions 2.1 and 2.5 we can read the value of $a_k(S/I)$'s for a monomial ideal $I$ of Borel type from its  irredundant  irreducible decomposition.

\vspace{2mm} {\bf \noindent Theorem  2.6.} {\it Let $I$ be a monomial
ideal of Borel type with a irredundant  irreducible decomposition
$I=\bigcap_{b\in B}{\mathfrak{m}}^{\bf{b}}.$  Then, for each $0\leq
k\leq n$, we have  $$a_k(S/I)=\max\{|{\bf{b}}|-n:  \
\mathrm{supp}({\bf{b}})=\{1,\cdots,n-k\}, {\bf{b}}\in B \}.$$ }
\noindent{\bf Proof.} We observe that
$I_1=I_0:x_{n_0}^{\infty}=\cap\{{\mathfrak{m}}^{\bf{b}}|\mathrm{supp}({\bf{b}})\subseteq
\{1,\cdots,n_1\},{\bf{b}}\in B\}$ and
$I_{l}=\cap\{{\mathfrak{m}}^{\bf{b}}:\mathrm{supp}({\bf{b}})\subseteq
\{1,\cdots,n_l\}, {\bf{b}}\in B\}, l=2,\cdots,r-1$ by induction.
Hence $\mathrm{Ass}(S/I)=\{(x_1,\cdots,x_{n_l}):l=0,\cdots,r-1\}$.
Fix $0\leq k\leq n$. If $k\notin \{n-n_0,\cdots,n-n_{r-1}\}$ then
the set $\{{\bf{b}}\in B: \mathrm{supp({\bf{b}})}=\{1,2,\cdots, n-k\} \}$ is empty, and so
$a_k(S/I)=\max\{|{\bf{b}}|-n: \
\mathrm{supp}({\bf{b}})=\{1, \cdots, n-k\} \}=-\infty$.
  If $k=n-n_l$ for some $l$, then $a_k(S/I)=s(J_l^{sat}/J_l)-k$. Let $B_1=\{{\bf{b}}\in B: \mathrm{supp}({\bf{b}})\subseteq\{1,\cdots,n_l\}\}$. Then
  $J_l=\cap \{\mathfrak{n}^{{\bf{b}}'}:{\bf{b}}\in B_1\}$,
  where $\mathfrak{n}^{{\bf{b}}'}=\mathfrak{m}^{{\bf{b}}}\cap k[x_1,\cdots,x_l],\ \forall {\bf{b}}\in B_1$.
   Since  $\mathrm{supp}({\bf{b}}')=\mathrm{supp}({\bf{b}})$, we have
    $a_k(S/I)=\max\{|{\bf{b}}|-n_l-k:\mathrm{supp}({\bf{b}})=\{1,\cdots,n_l\}, {\bf{b}}\in B_1\}=\max\{|{\bf{b}}|-n:  \
\mathrm{supp}({\bf{b}})=\{1,\cdots,n-k\}, {\bf{b}}\in B \} $ by
Proposition 2.5, as required. $\Box$

\vspace{2mm}

{\bf \noindent  Example 2.7.}  Let $I=(x^2,y^3)\cap
(x^4,y^4,z^3)=(x^4,x^2z^3,y^4,y^3z^3).$ Then $a_0(R/I)=8$ and
$a_1(R/I)=2$. Hence $\mathrm{reg}(R/I)=8$ and $a^*(R/I)=8$. \vspace{2mm}

The following corollary is immediate from Theorem 2.6.

\vspace{2mm} {\bf \noindent Corollary  2.8.} {\it  Let $I$ be a monomial
ideal of Borel type with a irredundant  irreducible decomposition
$I=q_1\cap q_2\cap \cdots \cap q_r.$  Then, for each $0\leq
k\leq n$, we have  $$a_k(S/I)=\max\{a_k(S/q_1),\cdots, a_k(S/q_r)\}.$$}

\vspace{2mm}

By Theorem 2.6, we can prove the main result of \cite{C1} easily.

\vspace{2mm}
{\bf \noindent Corollary  2.9.}\, (\cite[Corollary 8]{C1})  {\it Let $I$ be a monomial ideal of Borel type. Then $$\mathrm{reg}(I)=\min \{e\geq \mathrm{deg}(I)|I_{\geq e} \mbox{ is stable}\}.$$}
{\bf
\noindent Proof.} By \cite[Proposition 12]{ERT}, we only need to show if $e\geq \mathrm{reg}(I)$, then $I_{\geq e} \mbox{ is stable}$. We prove this result by induction on $r$, where $I=q_1\cap q_2\cap \cdots \cap q_r$ is  a irredundant  irreducible decomposition of $I$. It is clear in case that $r=1$. If $r>1$, we put $J=q_1\cap q_2\cap \cdots \cap q_{r-1}.$ Then ${\rm{reg}}(I)=\max \{{\rm{reg}}(J), {\rm{reg}}(q_r)\}$ by Theorem 2.6. It follows that for any $e\geq \mathrm{reg}(I)$, the ideals $J_{\geq e}$ and $(q_r)_{\geq e} \mbox{ are stable}$ by induction hypothesis and thus $I_{\geq e}=J_{\geq e}\cap (q_r)_{\geq e}$ is stable. as required. $\Box$ \vspace{2mm}

It is well-known that  $K\cap L$,  $K+L$  $KL$ are all of Borel type if  $K, L$ are two monomial ideals of Borel type.

\vspace{2mm} {\bf \noindent Corollary  2.10.} {\it Let $K,L$ be
two monomial ideals of Borel type. Set $I=K\cap L$. Then
$a_i(S/I)\leq \max\{a_i(S/K),a_i(S/L)\}, \ \forall 0\leq i\leq n$. In
particular}
$$\mathrm{reg}_t(S/I)\leq \max\{\mathrm{reg}_t(S/K),\mathrm{reg}_t(S/L)\}
\mbox{ and } a^*_t(S/I)\leq \max\{a^*_t(S/K),a^*_t(S/L)\}. $$
{\bf\noindent Proof.} Assume that $K=\bigcap_{b\in B_1}{\mathfrak{m}}^{\bf{b}},$ $L=\bigcap_{b\in B_2}{\mathfrak{m}}^{\bf{b}}$ and $I=\bigcap_{b\in C}{\mathfrak{m}}^{\bf{b}}$ are irredundant  irreducible decompositions of $K, L, I$ respectively. Note that $C\subseteq B_1\cap B_2$, we obtain that
$a_i(S/I)\leq \max\{a_i(S/K),a_i(S/L)\}, \ \forall 0\leq i\leq n$ immediately by Theorem 2.6. $\Box$

\vspace{2mm}

Let $K=(x^2, y^2, z^{10}), L=(x^4,y^4), I=K\cap L=L$. Then $a_0(S/I)=-\infty$ and $a_0(S/K)=11$. Hence, $a_i(S/I)=\max\{a_i(S/K),a_i(S/L)\}$ does not hold in general.

\vspace{2mm}

{\bf \noindent Corollary  2.11.} {\it Let $K,L$ be two monomial ideals of Borel type. Set $J=K+L$. Then}
$$  \mathrm{ (1)} \ a_t^*(S/J)=\max\{a_{t+1}^*(S/K),a_{t+1}^*(S/L)\}.  \,\, \quad \qquad \qquad \qquad \qquad \qquad \qquad \qquad $$
  $$ \mathrm{(2)} \  \mathrm{reg}_t(S/J)\leq \max\{\mathrm{reg}_t(S/K), \mathrm{reg}_t(S/L),\mathrm{reg}_{t+1}(S/K)-1,\mathrm{reg}_{t+1}(S/L)-1\}.$$
\vspace{2mm}
{\bf
\noindent Proof.} From the short exact sequence
$0\rightarrow S/(K\cap L)\rightarrow S/K\oplus S/L\rightarrow S/J\rightarrow 0$. We obtain the long exact sequence: $$\cdots \rightarrow H_{\mathfrak{m}}^i(S/K) \oplus H_{\mathfrak{m}}^i(S/L) \rightarrow H_{\mathfrak{m}}^i(S/J)\rightarrow  H_{\mathfrak{m}}^{i+1}(S/K\cap L) \rightarrow \cdots.$$ It follows that $a_i(S/J)\leq \max\{a_i(S/K),a_i(S/L) a_{i+1}(S/K),a_{i+1}(S/L) \}$ by Corollary 2.10. Now the reults follow directly.

\end{CJK*}

\end{document}